\documentclass[reqno]{amsproc}

\usepackage{amstext,amsmath,amssymb,amsfonts}
\usepackage[latin1]{inputenc}
\usepackage{epsfig}
\usepackage{hyperref}
\usepackage{color}

\usepackage{cancel}
\usepackage{latexsym}

\usepackage{enumitem}
\usepackage{bbm}
\usepackage{geometry}
\usepackage{amsthm}

\definecolor{myurlcolor}{rgb}{0,0,0.7}

\newtheorem{lemma}{Lemma}
\newtheorem{corollary}{Corollary}

\newtheorem{theorem}{Theorem}

\newcommand{\vspi}{\vspace{0.4cm}}

\newcommand{\bea}{\begin{eqnarray}}
\newcommand{\eea}{\end{eqnarray}}
\newcommand{\beq}{\begin{equation}}
\newcommand{\eeq}{\end{equation}}
\newcommand{\enn}{\nonumber \end{equation}}

 \newcommand{\cG}{\mathcal{G}}
 \newcommand{\cV}{\mathcal{V}}
 \newcommand{\cE}{\mathcal{E}}
 
 \newcommand{\bbN}{\mathbb{N}}
 \newcommand{\bbZ}{\mathbb{Z}}

\newcommand{\bin}[2]{\left (\begin{array}{c}#1\\#2\end{array} \right ) }

\def\eps{\epsilon}
\def\f{\frac}
\def\cN{{\mathcal N}}
\def\cC{{\mathcal C}}

\title[On terminal forms for topological polynomials 
for ribbon graphs]{On terminal forms \\ for topological polynomials 
for ribbon graphs: \\  The $N$-petal flower\footnote{Preprint: pi-mathphys-310, ICMPA-MPA/2012/35.}}

\author{ Remi C. Avohou}
\email{a) avohouremicocou@yahoo.fr}
\address{
International Chair in Mathematical Physics and Applications,
ICMPA-UNESCO Chair, 072BP50, Cotonou, Rep. of Benin
}

\author{Joseph Ben Geloun}
\email{b) jbengeloun@perimeterinstitute.ca}
\address{Perimeter Institute for Theoretical Physics, 31 Caroline
St, Waterloo, ON N2L 2Y5, Canada}
\address{International Chair in Mathematical Physics and Applications,
ICMPA-UNESCO Chair, 072BP50, Cotonou, Rep. of Benin}

\author{Etera R. Livine}
\email{c) etera.livine@ens-lyon.fr}
\address{Laboratoire de Physique, ENS Lyon, CNRS-UMR 5672, 46 All\'ee d'Italie, Lyon 69007, France}
\address{Perimeter Institute for Theoretical Physics, 31 Caroline
St, Waterloo, ON N2L 2Y5, Canada}

\begin{document}

\maketitle 

\begin{abstract}
The Bollobas-Riordan polynomial [Math. Ann. 323, 81 (2002)] 
extends the Tutte polynomial and its contraction/deletion rule
for ordinary graphs to ribbon graphs. 
Given a ribbon graph $\cG$, the related polynomial should
be computable from the knowledge of the 
terminal forms of $\cG$ namely specific
induced graphs for which the contraction/deletion procedure becomes more involved. 
We consider some classes of terminal forms 
as rosette ribbon graphs with $N\ge 1$ petals 
and solve their associate Bollobas-Riordan polynomial. 
This work therefore enlarges the list of terminal 
forms for ribbon graphs  
for which the Bollobas-Riordan polynomial
could be directly deduced.   \\ 

 \end{abstract}

\tableofcontents

\section{Introduction: Background and motivations}

Bollobas and Riordan in \cite{bollo,bollo2} introduced a polynomial
for  ribbon graphs or graphs on surfaces. 
Let us review few ingredients necessary to their analysis. 

A ribbon graph $\cG$ is  a (not necessarily orientable) surface with boundary represented as the union of 
two sets of closed topological discs called vertices $\cV$ and edges $\cE.$ These sets satisfy the following:
 Vertices and edges intersect by  disjoint line segments;
 each such line segment lies on the boundary of precisely one vertex and one edge, every edge contains exactly two such line segments.

An edge $e$ of a graph $\cG$ can have specific properties: 
 $e$ is called a self-loop in $\cG$ if the two ends of $e$ are adjacent to the same vertex $v$ of $\cG$;
 $e$ is called a bridge in $\cG$ if its removal disconnects a component of $\cG$;
 $e$ is called an ordinary or regular edge of $\cG$ if it is neither a bridge nor a self-loop.
A graph which does not contain any regular edge 
shall be called a terminal form. 

\begin{figure}
 \centering
     \begin{minipage}[t]{.8\textwidth}
      \centering
\includegraphics[angle=0, width=5cm, height=2cm]{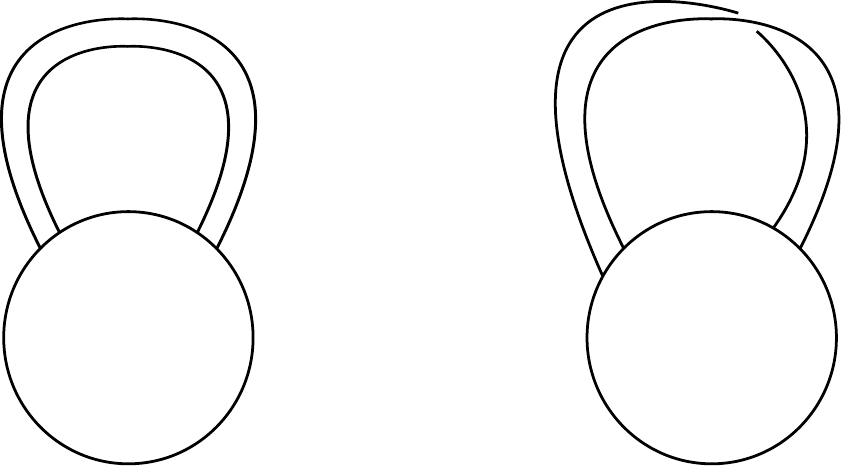}
\caption{ {\small Untwisted (left) and twisted (right) trivial self-loop. }}
\label{fig:twist}
\end{minipage}
\end{figure}

Focusing on particular self-loops, one has \cite{bollo}:
 A self-loop $e$ at some vertex $v$ is called trivial
if there is no loop $f$ at $v$ such that the ends of $e$ and $f$ 
alternate in the cyclic order at $v$. 
A loop $e$ at $v$ is called twisted if $v \cup e$ 
forms a M\"obius band; if $v \cup e$ forms an annulus 
$e$ is called a untwisted self-loop (see Figure \ref{fig:twist}). 

The notion of contraction and deletion of an edge \cite{tutte,bollo} (see also \cite{riv,krf} for interesting connections with quantum 
field theory) is now recalled. 
Let $\cG$ be a graph and $e$ one of its edges. 
 We call $\cG-e$ the graph obtained from 
$\cG$ by removing $e$.
 If $e$ is not a self-loop, the  graph $\cG/e$  obtained by contracting $e$ is defined from $\cG$ by deleting $e$ and 
identifying its end vertices into a new vertex;
If $e$ is a self-loop, $\cG/e$ is by definition the same as $\cG-e$.

One notices that after a contraction-deletion sequence of all ordinary edges of a given graph 
the end result is necessarily given by a collection 
of graphs composed by  bridges and/or self-loops, hence 
a terminal form. 

The Bollobas-Riordan (BR) topological polynomial for 
ribbon graph is given by \cite{bollo} 
\beq
R_{\cG}(X,Y,Z) =  \sum_{A \subset \cG} 
(X-1)^{r(\cG) - r(A)} 
Y^{n(A)} Z^{k(A)- bc(A)+ n(A)}\,,
\eeq
where the sum is over all spanning subgraphs $A$ of $\cG$, 
and using standard parameters for graph \cite{tutte,bollo}
$v(A)= v(\cG)$ is the number of vertices of $A$, 
$E(A)$ is the number of edges of $A$, 
$k(A)$ is the number of connected components
 of $A$, $r(A)$ is the rank of $A$ and is given by 
$r(A) = v(\cG) -  k(A)$, $n(A) = E(A) - r(A)$
is the nullity of $A$ (or first Betti number). 
In addition, $bc(A)$ is the number of components of the
boundary of $A$ when $A$ is regarded as a geometric ribbon graph
\cite{bollo}. We simply call it number of faces $bc(A) = F(A) $ in the following.  

Note that in \cite{bollo} there is an extra variable
in $W$ which takes into account the orientability 
of the subgraph when seen as a surface. We simply put this variable to 1 in the present analysis. 
Our result should find an extension for general $W$. 

The BR polynomial $R_{\cG}$ is called topological because it satisfies the following contraction and deletion rules. For an ordinary edge $e$, 
we have
\beq
R_{\cG }  = R_{\cG - e} + R_{\cG /e}\,, 
\eeq
for every bridge $e$ of $\cG$, 
\beq
R_{\cG}=X \, R_{\cG/e}\,,
\label{cbrid}
\eeq 
for a trivial untwisted self-loop, 
\beq
R_{\cG}=(1+Y) \,R_{\cG-e}\,,
\label{csel}
\eeq
and for a trivial twisted self-loop, the following holds
\beq
R_{\cG}=(1+YZ) \,R_{\cG-e}\,.
\label{ctsel}
\eeq
The relations  \eqref{cbrid}-\eqref{ctsel} 
are useful for the evaluation of the BR polynomial of
a graph $\cG$ from its terminal forms. For instance, if at the end
of a contraction and deletion sequence of all regular edges, a ribbon 
graph $\cG$ yields some disconnected graph family  $\{\cG_i\}$ with each
 $m_i$ bridges, $p_i$ trivial untwisted and $q_i$ trivial twisted self-loops then 
the BR polynomial associated with such a graph $\cG$
will be simply a summation of the contributions 
\beq\label{contribself}
X^{m_i} (1+Y)^{p_i} (1+YZ)^{q_i}\,. 
\eeq

However, the above listed terminal forms for ribbon graphs 
are far to be exhaustive. 
It noteworthy that the Tutte polynomial 
for a graph $\cG$ can be always  evaluated 
from contraction and deletion moves applied to only regular edges of $\cG$ yielding computable terminal
forms\footnote{The resulting disconnected family $\{\cG_i\}$ of terminal forms have a Tutte polynomial 
which is directly computable in terms of $X^{m_i} (1+Y)^{p_i}$
where $m_i$ is the number of bridges and $p_i$ 
the number of self-loops.}. In contrast, in a generic situation,
after the full contraction and deletion sequence of all regular edges, the BR polynomial of a ribbon graph may be not directly evaluated because all terminal forms including a subgraph of the form of a rosette graph (single vertex ribbon graph) have been not yet solved. 
In last resort, the contribution of these terminal forms of the ribbon graph can be only computed through the summation over subgraphs. 
 A natural question follows: ``Is it possible to enlarge the space of computable `initial conditions' by providing an explicit BR polynomial expression of the most general rosette graph?''. This question is more intricate than one might
think because the BR polynomial is more than
 a simple topological invariant of a ribbon graph 
 (this is also the case for the Tutte polynomial for graphs). 
One points out also that 
even if a one-vertex ribbon graph may be mapped to a
much simpler specific ribbon graph  
(via the so-called chord diagrams crucial in the proof of the universality property of the BR polynomial \cite{bollo2, bollo} and useful in 
the other context of Vassiliev invariants, see for instance \cite{birman}) there is no clear way to extract from these simplest configurations 
the BR associated to the anterior graph itself. 
In this work, the above question
finds a partial but positive answer by solving a less stronger problem 
for specific classes of rosette graphs.

In this paper, we study some families of single vertex ribbon graphs (parametrized by 
their number of edges supplemented by other features) with self-loops  
the members of which should be considered each as a terminal form
different from the aforementioned (trivial twisted and untwisted self-loops). 
Thus, we aim at completing the contribution \eqref{contribself}
if some specific members of these families occur as part of a graph 
inferred by contraction and deletion sequence of all regular
edges.  Interestingly, we find that the face counting in subgraphs 
of these families of graphs involve some number of specific compositions
(ordered partitions). This certainly foresees rich links between number theory and topological polynomials. 
   
\begin{figure}
 \centering
     \begin{minipage}[t]{.8\textwidth}
      \centering
\includegraphics[angle=0, width=5cm, height=2cm]{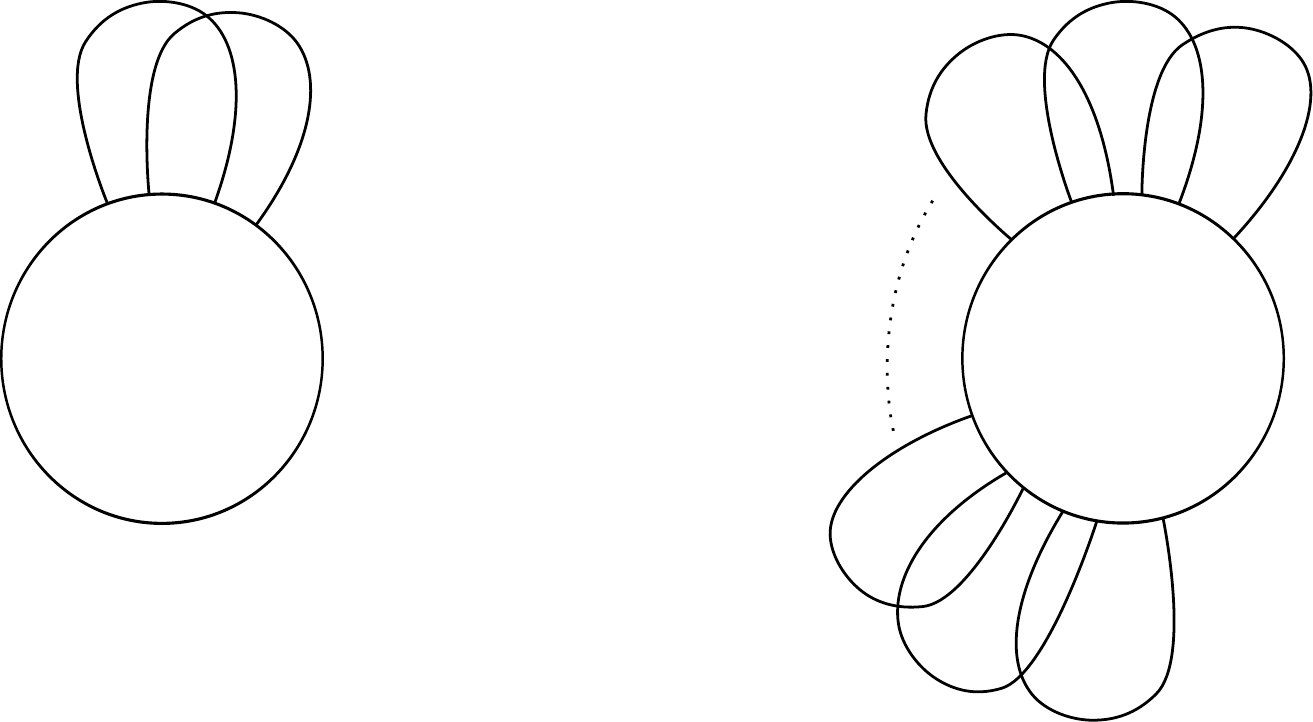}
\caption{ {\small A graph $\cG_2$ with two non trivial self-loops
and 
the resulting $N$-petal flower graph $\cG_N$. }}
\label{fig:nflow}
\end{minipage}
\put(-228,0){$\cG_2$}
\put(-100,0){$\cG_N$}
\end{figure}

\begin{figure}
 \centering
     \begin{minipage}[t]{.8\textwidth}
      \centering
\includegraphics[angle=0, width=12cm, height=3cm]{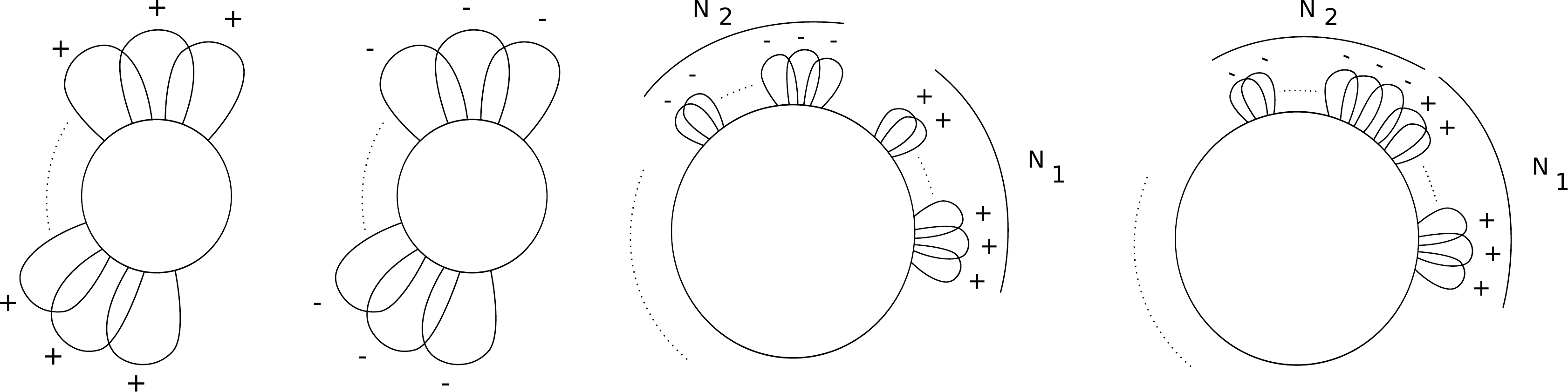}
\vspace{0.2cm}
\caption{ {\small  The $N$-twisted-petal flower $\cG_{N}^+$,  the $N$-untwisted-petal flower $\cG_{N}^-$ and 
the two generalized $\{(N_1,+),(N_2,-),\dots,(N_q,s_q)\}$-petal
flowers $\cG_1$ and $\cG_2$ with seperate and merged
sectors, respectively.}}
\label{fig:Nsflow}
\end{minipage}
\put(-325,-10){A: $\cG^+_{N}$}
\put(-255,-10){B: $\cG^-_{N}$}
\put(-175,-10){C: $\cG_1$}
\put(-77,-10){D: $\cG_2$}
\end{figure}
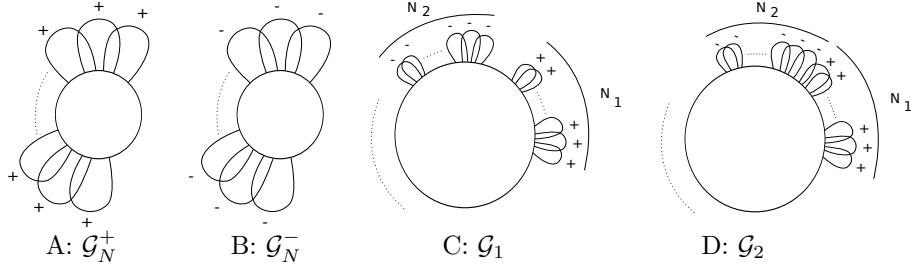

Consider $\cG_2$ in Figure \ref{fig:nflow} (we use there and 
in Figure \ref{fig:Nsflow} simplified diagrammatics where each edge should be viewed as a ribbon). This is a rosette graph which has two non trivial self-loops. The same diagram extends to a special rosette graph 
that we will be referring to $N$-petal flower (see $\cG_N$ in Figure \ref{fig:nflow}). 
Note that the edges are intertwined in a specific way. Given the fact that each petal  
can be twisted or not, we have a variety of different flowers. 
Ultimately, we will be interested in the more elaborate
$\{(N_1,s_1),(N_2,s_2),\dots, (N_q,s_q)\}$-petal flower 
where each sector $(N_l,s_l)$ refers to a number $N_l$
of petals which can be twisted $(s_l=+)$ or not $(s_l=-)$. 
Such a flower can be a separate collection of $(N_l,s_l)$ sectors  
(Figure \ref{fig:Nsflow} C)  
or a unique collection of  $(N_l,s_l)$ sectors
merged to their neighbor(s) (Figure \ref{fig:Nsflow} D) .

\vspi
We aim at computing the BR polynomial for 

(1) the  $N$-untwisted-petal
flower (see Figure \ref{fig:Nsflow} B), with $N \in \mathbb{N},$ 
and $N \geq 1$;

(2) the $N$-twisted-petal
flower (see Figure \ref{fig:Nsflow} A) with  $N \in \mathbb{N}$
and $N \geq 1$;

(3) the generalized $\{(N_1,s_1),(N_2,s_2),\dots, (N_q,s_q)\}$-petal
flower (see Figure \ref{fig:Nsflow} C) with  $N_l\in \mathbb{N}$, $N \geq 1$, $1 \le l \le q$. Note that each sector $(N_l,s_l)$ is not connected to any 
other sector.

This can be achieved after finding an explicit formula 
for their number of faces. Hence, we compute in a more general 
setting and useful in any situation,

(0) the number of faces 
of the generalized 
$\{(N_1,s_1),(N_2,s_2),\dots,(N_q,s_q)\}$-petal flower
 given by Figure \ref{fig:Nsflow} D.

We call trivial a $N$-petal flower at some vertex $v$, 
if the only possible ends of any loop which alternate in the cyclic order
at $v$ with the ends of the flower edges  belongs to the flower itself.

\section{Main results}

In this paper, we prove the following statements
which are our main results: 

\begin{theorem}[Number of faces of the $N$-petal flower]
\label{theo:face}
Given a $N$-petal flower with twisted and untwisted  
petals. The number of faces of this graph is 
either 1 or 2. 
\end{theorem}

\begin{theorem}[BR polynomial for the $N$-(un)twisted petal flower]
\label{theo:brNu}
Given a $N$-untwisted petal flower $\cG_N$, the BR polynomial associated
with  $\cG_N$ is given by 
\beq\label{eq:unflow}
R_{\cG_N}(Y,Z)
\,=\,
1+
\sum_{n=1}^{N-1} Y^n\, \sum_{P=1}^n \bin{N-n+1}{P}
\sum_{I=\eps(n)}^{P} Z^{n-I} \cC_n(P,I)
+Y^NZ^{N-\eps(N)}\,,
\eeq
where $\eps(q) = (1-(-1)^q)/2$, for any $q\in \mathbb{N}$,
 and 
$\cC_{n}(P,I)$ is the number of compositions  (a.k.a. ordered partitions) 
of the integer $n$ in $P$ integers among which $I$ odd integers.

Given a $N$-twisted petal flower $\cG^t_N$, the BR polynomial associated
with  $\cG^t_N$ is given by 
\beq\label{eq:twflow}
R_{\cG^t_N}(Y,Z)
\,=\,
1+
\sum_{n=1}^{N-1} Y^n\, \sum_{P=1}^n \bin{N-n+1}{P}
\sum_{I=0}^{P} Z^{n-I} \cC^t_n(P,I)
+Y^NZ^{N-\eps_{(3)}(N)}\,,
\eeq
where $\eps_{(3)}(q) = 1$ if $q \in 3\bbN +2$, otherwise $\eps_{(3)}(q)=0$,  
$\cC^t_{n}(P,I)$ is the number of compositions of the integer $n$ in $P$ integers among which $I$ integer belonging to $3\bbN +2$.
\end{theorem}

Setting $Z=1$ in the polynomials \eqref{eq:unflow} 
and \eqref{eq:twflow} one simply recovers the Tutte polynomial $T_{\cG_N}$ for a simple graph with $N$ self-loops namely
\beq
R_{\cG_N}(Y,Z=1)=R_{\cG^t_N}(Y,Z=1)= (1+Y)^N= T_{\cG_N}(Y)\,.
\eeq
Hence,  at  $Z=1$, the coefficient of $Y^n$ in these equations can be
simply regarded as a peculiar decomposition of the binomial 
coefficient $(^N_n)$ in terms of number of particular compositions.

\begin{corollary}[Deleting a $N$-petal flower]
\label{coro:delflo}

Let $\cG$ be a graph with a trivial $N$-petal flower subgraph 
$G_{N}$ which could be twisted or not. Deleting all petals of $G_{N}$, we define the resulting graph $\cG-G_N$
and we have 
\bea
R_{\cG} = R_{G_{N}} R_{\cG - G_N}\,,
\eea
where $R_{G_N}$ is obtained by \eqref{eq:unflow}
if $G_N$ is $N$-untwisted and given by \eqref{eq:twflow} 
if $G_N$ is $N$-twisted. 

\end{corollary}

\begin{corollary}[BR polynomial for generalized terminal forms]
\label{coro:terminal}

Let $\cG$ be a graph made with 
$m$ bridges, $p$ trivial untwisted self-loops, 
$q$ trivial  self-loops, a finite family $\{\cG_{N_i}\}_{i\in I}$
of trivial $N_i$-untwisted-petal flowers and a finite family 
$\{\cG^t_{N_j}\}_{j\in J}$ of trivial $N_j$-twisted-petal
flowers. The BR polynomial for $\cG$ is given by 
\beq
X^m (1+Y)^{p}(1+YZ)^{q} 
\bigg[\prod_{i\in I} R_{\cG_{N_i}} \bigg]
\bigg[\prod_{j\in J} R_{\cG^t_{N_j}} \bigg]\,.
\eeq

\end{corollary}

\section{Number of faces of a generalized $N$-petal flower}
\label{sect:gene}

This section is mainly devoted to the proof of Theorem 
\ref{theo:face}. Furthermore, we investigate  
useful consequences of this result. 

Let us emphasize first that it may exist another proof
of this statement using chord diagrams $D_{ij}$ used in  \cite{bollo,bollo2}. For flower untwisted and twisted petals,
this may be quickly achievable. However recasting the generalized
situation in terms of these ``canonical'' chord diagrams (taking into
account the orientations induced by twistings) shall need
a non trivial algorithm and so the proof of Theorem \ref{theo:face}
should be in any way non trivial. We will use another method
which is  itself interesting.

Consider a generalized 
$\{(N_1,s_1),(N_2,s_2),\dots, (N_q,s_q)\}$-petal flower given 
Figure \ref{fig:Nsflow} D 
 with the specific feature that each sector $(N_l,s_l)$ is  
connected to its neighbor sector(s). 
We simply refer such a rosette, in this section, to as $N$-petal
flower. 

\begin{figure}
 \centering 
     \begin{minipage}[t]{.8\textwidth}
      \centering
\includegraphics[angle=0, width=3cm, height=2.5cm]{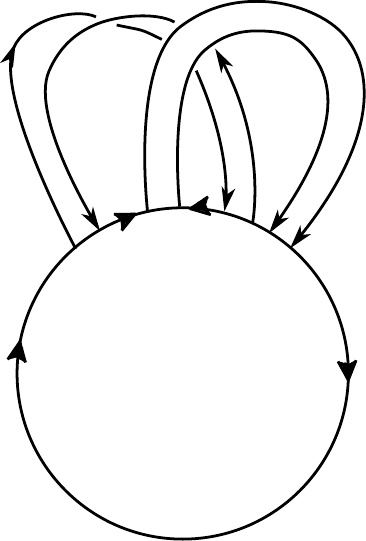}
\vspace{0.2cm}
\caption{ {\small  An example of  flower with one oriented face. }}
\label{fig:orient}
\end{minipage}
\end{figure}

Counting the number of faces of at some fixed 
 $N$ number of edges of the $N$-untwisted or 
twisted flower can be simply achieved by induction. 
However, for a general $N$-petal flower, the number of faces 
becomes intricate. 
A way to overcome this issue is to introduce another ingredient
on the graph which is the notion of orientation of each face. An orientation
is simply denoted as an arrow on the face, see Figure \ref{fig:orient}. 
This corresponds to an orientation (in the geometric sense)
of the boundary of the ribbon graph 
when the graph is viewed as a geometric ribbon. 
Note that  this type of orientation should be related to the edge orientation in the sense defined in \cite{bollo}. In any case, a face orientation induces an edge orientation,  that is an orientation 
of its side segments.  We say that a graph has a face
orientation if to all of its faces we assign an arrow. 

We emphasize that one can identify for a $N$-petal 
flower an initial and last petal given a cyclic order
on the vertex. 

Given a $N$-petal flower equipped with a face orientation, 
the number $F$ of faces for such a graph can be obtained
 by recurrence on the number of petals in the flower and
the orientations of the two sides of the last end of the last petal
(or the first, without loss of generality). 

\proof[Proof of Theorem \ref{theo:face}]
Theorem \ref{theo:face} claims that
the number of faces of the $N$-petal flower is 1 or 2.
In order to prove this by recurrence, we adopt the following strategy: 
at the order $N$, we add another petal with a given edge orientation and observe the change in the number of faces of the resulting 
flower. Note that we also need to encode the change in the
orientation induced by the additional edge. 

Let us quickly see how the previous statement translates
for the lowest orders for $N>0$. 
Consider $N=1$, this is a 1-petal flower so that either 
$F=1$ (twisted petal) or $F=2$ (untwisted petal). 

\begin{figure}[h]
 \centering
     \begin{minipage}[t]{.8\textwidth}
      \centering
\includegraphics[angle=0, width=10cm, height=3cm]{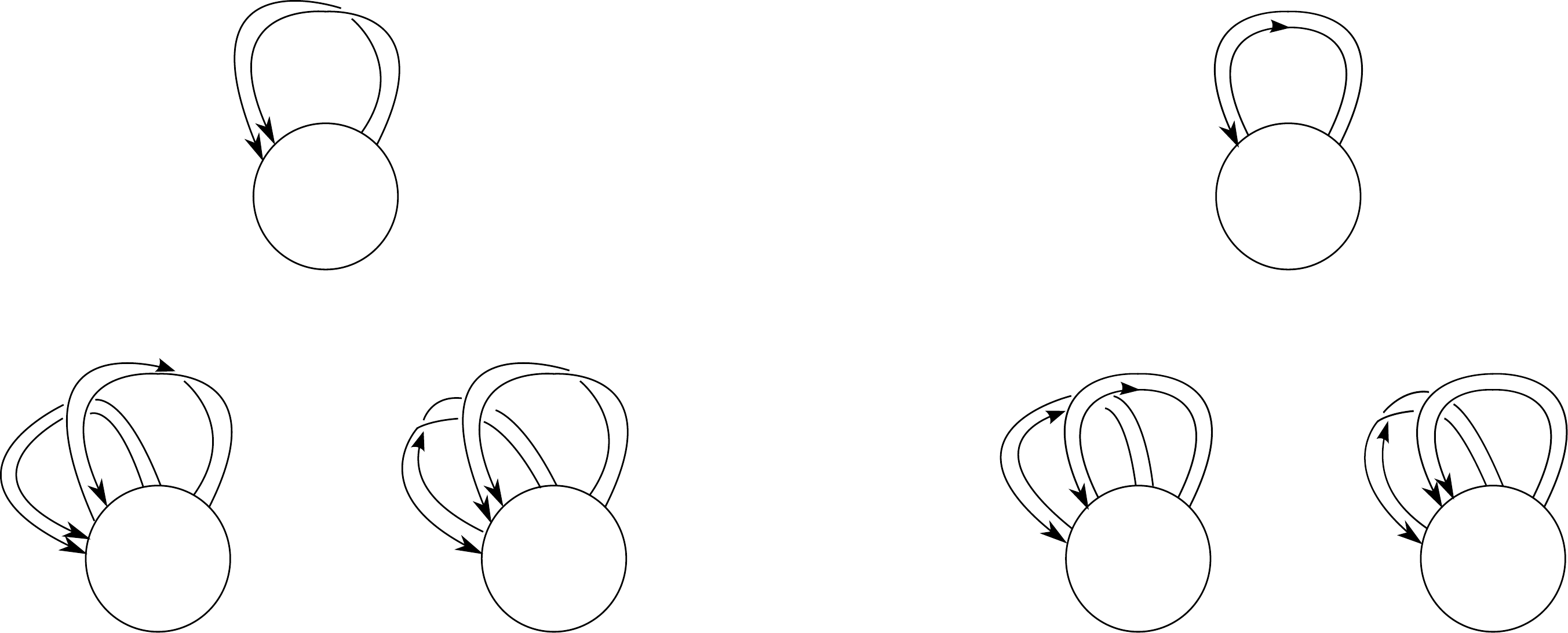}
\vspace{0.2cm}
\caption{ {\small  $(N=1)$-twisted petal flower (A) 
with one oriented face and insertion of an untwisted (AB) and twisted (AA)
petal;  $(N=1)$-untwisted petal flower (B) 
with oriented faces and insertion of an untwisted (BB) and twisted (BA)
petal.}}
\label{fig:G1twist}
\end{minipage}
\put(-257,37){$A$}
\put(-293,-10){$AB$}
\put(-220,-10){$AA$}
\put(-84,37){$B$}
\put(-115,-10){$BB$}
\put(-50,-10){$BA$}
\end{figure}
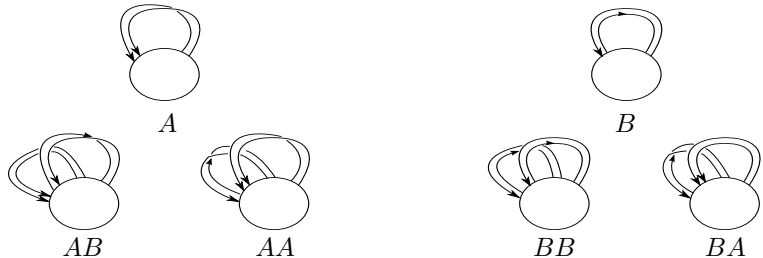

(A) For the twisted petal, $F=1$: assume an initial (and unique) 
orientation called $(s=+)$ of the face given by Figure \ref{fig:G1twist} A. 

(AB) Adding  an untwisted petal gives Figure \ref{fig:G1twist} $AB$
with $F=1$ and observe that the last end of the edge in $AB$
possesses the same orientation $(s=+)$;

(AA) Adding a twisted petal yields Figure \ref{fig:G1twist} $AA$
with $F = 2$. Each side of the last edge does not belong 
to the same face, so that their orientation $s$ is arbitrary. 
In such a situation, we do not need to report the edge orientation.

(B) For the untwisted petal, $F=2$: the orientations
of the faces are arbitrary. We choose for instance those given by Figure \ref{fig:G1twist} B. One checks that the following is independent
of that initial face orientations. 

(BB) Adding  an untwisted petal leads to Figure \ref{fig:G1twist} $BB$
with $F=1$ and observe that the last end of the edge in $BB$
possesses an orientation that we will be referring to $(s=-)$
as opposed to the $(+)$ of Figure \ref{fig:G1twist} $A$;

(BA) Adding a twisted petal gives Figure \ref{fig:G1twist} $BA$
with $F = 1$. Each side of the last edge belongs 
to the same face with orientation $(s=-)$.

Adding more petals to these configurations, one rapidly 
finds the recurrence hypothesis as follows: given a flower
with $N$ petals with $F=1,2$ face(s)  and $s=\pm$ 
edge orientation of the last petal, adding 
a new twisted petal $(\text{t})$ or an untwisted petal $(\text{unt})$
yields the table:
\bea
(F=1,s=+)\quad  \stackrel{\text{t}}{\longrightarrow} \quad (F=2,s)
\,,&&
(F=1,s=+)\quad  \stackrel{\text{unt}}{\longrightarrow} \quad (F=1,s=+)\,,
\crcr
(F=1,s=-)\quad  \stackrel{\text{t}}{\longrightarrow} \quad (F=1,s=+)
\,,&&
(F=1,s=-)\quad  \stackrel{\text{unt}}{\longrightarrow} \quad (F=2,s)\,,
\crcr
(F=2,s)\quad  \stackrel{\text{t}}{\longrightarrow} \quad (F=1,s=-1)
\,,&& 
(F=2,s)\quad  \stackrel{\text{unt}}{\longrightarrow} \quad (F=1,s=-1)\,.
\label{eq:hypo}
\eea
Note that our previous test on $N=1,2$ yields the first 
and last line in \eqref{eq:hypo}.  
For the middle line, one has to add a petal to the case $N=2$,
hence going to $N=3$,
in order to obtain such an occurrence. For instance, starting from 
$BB$ or $BA$ and adding another petal, one recovers the second
line relations. By convention, we equip the simple vertex 
graph (a disc) with the data $(F=1,s=-1)$ so that 

- adding a twisted petal 
yields $(F=1,s=+)$ and thereby rejoining the configuration 
of Figure \ref{fig:G1twist} A

- and by adding a untwisted petal one gets  $(F=2,s)$
as it should be for this configuration, see Figure \ref{fig:G1twist} B. 

We assume that these relations holds at order $N$. 
Let us prove them for the  $(N+1)$-petal flower. 
For this purpose,  a case by case study is required. 

\begin{figure}[h]
 \centering
     \begin{minipage}[t]{.8\textwidth}
      \centering
\includegraphics[angle=0, width=8cm, height=3cm]{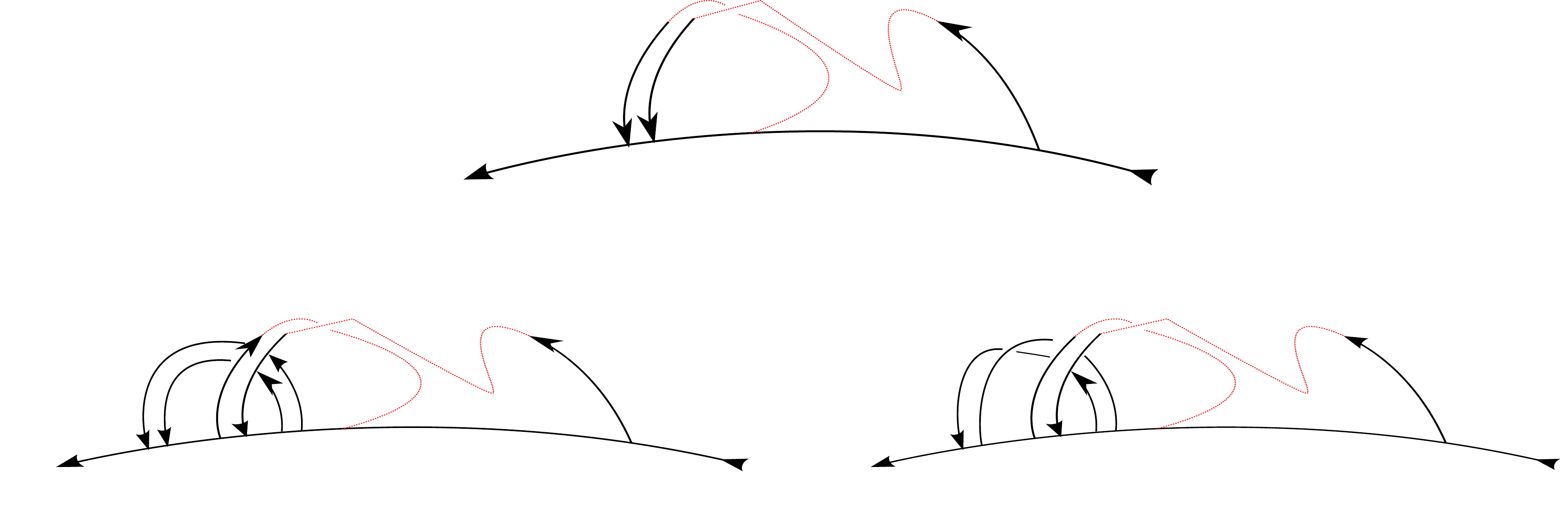}
\vspace{0.2cm}
\caption{ {\small  
A generalized  flower $(F=1,s=+)$ 
and two possible insertions of an untwisted (1A) and a twisted (1B)
petal yielding $(F=1,s=+)$ and $(F=2,s)$, respectively.}}
\label{fig:GNtwist}
\end{minipage}
\put(-277,65){$(F=1,s=+)$}
\put(-235,-10){$1A$}
\put(-115,-10){$1B$}
\end{figure}

(1) Let us start by a configuration with $(F=1,s=+)$. 
Note that the unique face is necessarily of the form given 
by  Figure \ref{fig:GNtwist} (where the particular paths in red
in the flower are irrelevant for the analysis but
only are relevant the different connections between
the last edge and the red paths). 

(1A) Adding to this flower an untwisted petal, 
one gets the configuration of Figure \ref{fig:GNtwist} 1A, 
so that $(F=1,s=+)$. 

(1B) Adding to the flower a twisted petal, 
 the configuration of Figure \ref{fig:GNtwist} 1B is obtained 
so that $(F=2,s)$.

\begin{figure}[h]
 \centering
     \begin{minipage}[t]{.8\textwidth}
      \centering
\includegraphics[angle=0, width=8cm, height=3cm]{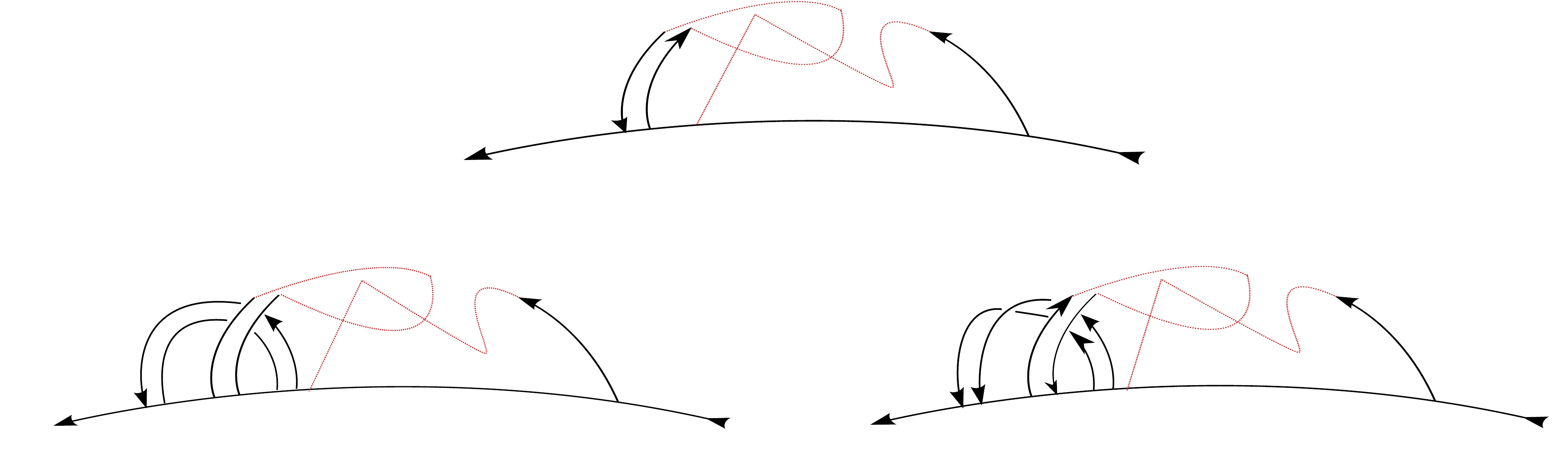}
\vspace{0.2cm}
\caption{ {\small  
A generalized  flower $(F=1,s=-)$ 
and two possible insertions of an untwisted (2A) and a twisted (2B)
petal yielding $(F=2,s)$ and $(F=1,s=+)$, respectively.}}
\label{fig:GNtwist2}
\end{minipage}
\put(-277,65){$(F=1,s=-)$}
\put(-235,-10){$2A$}
\put(-115,-10){$2B$}
\end{figure}

 (2) We pursue with the configuration $(F=1,s=-)$
described by the unique face which is necessarily 
of the form given by  Figure \ref{fig:GNtwist2}. 

(2A) Adding to this flower an untwisted petal, 
one gets  Figure \ref{fig:GNtwist} 2A, 
so that $(F=2,s)$.

(2B) Meanwhile, adding to the flower a twisted petal, 
one ends up with  Figure \ref{fig:GNtwist} 2B, 
so that $(F=1,s=+)$.

\begin{figure}[h]
 \centering
     \begin{minipage}[t]{.8\textwidth}
      \centering
\includegraphics[angle=0, width=8cm, height=3cm]{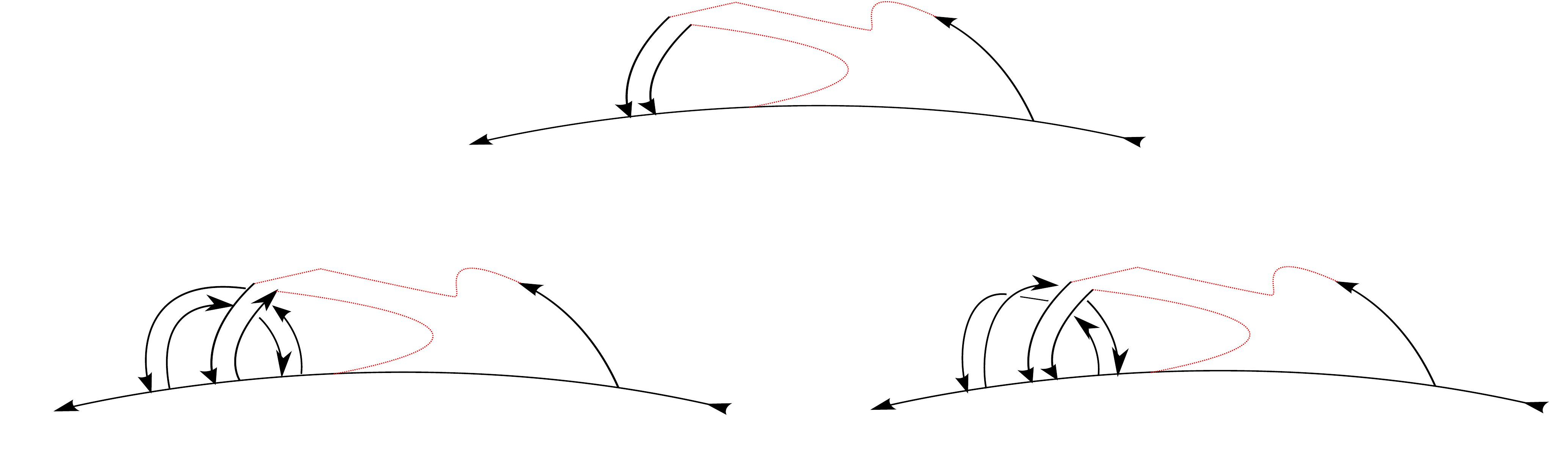}
\vspace{0.2cm}
\caption{ {\small  
A generalized  flower $(F=2,s)$ 
and two possible insertions of an untwisted (3A) and a twisted (3B)
petal yielding the same configuration $(F=1,s=-)$.}}
\label{fig:GNtwist3}
\end{minipage}
\put(-277,65){$(F=2,s)$}
\put(-235,-10){$3A$}
\put(-115,-10){$3B$}
\end{figure}

 (3) Finally, we study the configuration such that $(F=2,s)$
described by two faces with arbitrary orientations 
of the form given by  Figure \ref{fig:GNtwist3}. 

(3A) and (3B) Adding to this flower either an untwisted
or a twisted petal, the result is given in 
 Figure \ref{fig:GNtwist} 3A or 3B, respectively.
In any situation, one finds $(F=1,s=-)$. 

Then all relations \eqref{eq:hypo} are satisfied. 
Hence starting with any configuration of the $N$-petal flower 
and adding a new petal, either yields $F=1$ or $F=2$.
This achieves the proof of Theorem \ref{theo:face}. 

\qed 

We have thus obtained the number of faces of the $N$-petal flower
with arbitrary number of (twisted or untwisted) petals and without the explicit dependence
on the face orientations. Orientations have been used as an artifact of
the procedure. However, it is not clear that the number of faces is
dependent or not of the cyclic order of the vertex and the fact that 
we can distinguish two special edges, the first and the final, according to that cyclic order. With commonsensical arguments, one may claim that above 
number of faces should be the same if we interchange
the role of these two edges and reverse the cyclic order of the vertex. 
 We will come back on this point later and prove that
this is indeed the case.  

Theorem \ref{theo:face} is only useful if 
an explicit formula for the number of faces of
a $N$-petal flower is affordable. 
It is remarkable that we can map the each configuration 
on a $\bbZ_3$ generator and \eqref{eq:hypo} can be simply 
encoded in terms of a rule for these generators. We assign

-  $(F=1,s=+) \longrightarrow   2 \mod 3$

-  $(F=1,s=-) \longrightarrow   1 \mod 3$

-  $(F=2,s) \longrightarrow   0 \mod 3$

Notice that, according to our convention, 
the bare vertex gets mapped on $1$.
Defining for any petal $e$ the symbol 
$\varepsilon_e=-1$ if the petal is untwisted or 
$\varepsilon_e= +1$ if the petal is twisted, 
the above \eqref{eq:hypo} rules simply translate
as
\bea
x' = x \varepsilon_e + 1 \mod 3\,,
\eea
where $x$ is $\bbZ_3$ generator corresponding to 
the couple $(F,s)$ as defined above. 
For instance, for $x=2$, equivalently $(F=1,s=+)$,
calculating $ 2 (+1) + 1 = 3 =0$ equivalently 
describes the addition of a twisted petal to a configuration 
$(F=1,s=+)$ which yields $(F=2,s)$.

Consider now a  
$N$-petal flower as the result of
 branching $N \ge 1$ petals on an initial 
bare vertex. The bare vertex provides us
with an initial condition $x_0=1\in \bbZ_3$. 
Inserting a first petal $e_1$ we obtain the class
$x_1\in \bbZ_3$, $x_1 = x_0\varepsilon_{1}+1$, where, for simplicity, 
we henceforth denote $ \varepsilon_{e_i}=  \varepsilon_{i}$.
Then we iterate the procedure by inserting more petals such that,
at the end, the number of faces of the flower is directly obtained  
after evaluating a nested product 
\bea
(F,s) &=& [[\dots [[[\varepsilon_{1}+1]\varepsilon_2 + 1]\varepsilon_{3}+1]\dots ] \varepsilon_{N-1}+1]\varepsilon_{N} + 1
\mod 3\crcr
&=& 
1+ \sum_{l=0}^{N-1} \prod_{k=0}^l\varepsilon_{N-k}  \mod 3\,.
\label{eq:fac}
\eea 
Consider an ordered sequence $e_1, \dots, e_l,\dots e_N$ of self-loops
forming the $N$-petal flower, the subset 
$E_{l} = \{e_N,\dots,e_l\}$ of the final self-loops (starting from 
a final edge $e_N$ and counting in some cyclic order at the vertex up to $e_1$) and the subset $E_{\text{unt};l}\subset E_l$ of all  its untwisted self-loops,
then the class $(F,s)$ \eqref{eq:fac} can be rewritten as
\bea
(F,s) = 1+ \sum_{l=0}^{N-1} (-1)^{|E_{\text{unt},N-l}|} \mod 3\,.
\label{eq:fs}
\eea
Now we come back on a previous remark  and consider the reverse construction of the flower. We start by inserting the last edge $e_N$ 
then $e_{N-1}$ and so on up to $e_1$. Thus the class that one obtains
is  
\bea
(F',s')&=&  [[\dots [[[\varepsilon_{N}+1]\varepsilon_{N-1} + 1]\varepsilon_{N-2}+1]\dots ] \varepsilon_{2}+1]\varepsilon_{1} + 1
\mod 3\crcr
&=& 1+
\sum_{l=1}^{N}\prod_{k=1}^l\varepsilon_{k} \mod 3\,.
\label{eq:fac2} 
\eea
It is direct to recover 
\bea
(F',s') = (-1)^{|E_{\text{unt},1}|} (F,s) \mod 3\,. 
\eea
Hence, $(F,s)  = 0  \Leftrightarrow (F',s') = 0$ from which one infers
that $F=2=F'$, and otherwise necessarily $F=1=F'$. 

Let us restrict the formula \eqref{eq:fs} for particular flowers. 
Assuming that the flower has only untwisted petals, then, 
according to our rules, by adding successively petals only
the sequence $\dots \to 0\to 1\to 0 \to 1\dots $ is possible. 
Moreover,
\bea
(F,s) &=& 1+ \sum_{l=0}^{N-1} (-1)^{|E_{N-l}|} \mod 3
=  1+ \sum_{l=0}^{N-1} (-1)^{l+1} \mod 3 \crcr
& =&  \frac12 (1+(-1)^{N}) \mod 3\,.
\eea
Thus, if $N$ is even $(F,s)=1 \mod 3$ and 
if $N$ is odd $(F,s)=0 \mod 3$. 
In any situation, we directly identify $F=1+\eps(N)$,
where $\eps(N)=(1-(-1)^N)/2$.  

Assuming that the flower has only twisted petals, then 
for $N \geq 1$, 
\beq
(F,s) = 1+ \sum_{l=0}^{N-1} (-1)^{0}\mod 3=( 1+  N) \mod 3\,.
\eeq
In this case, we simply write $ F = 1+\eps_{(3)}(N)$.

 Note that formula \eqref{eq:fs} 
determines the number of faces of more general classes of $N$-petal flowers than the two simplest situations discussed above. For instance, whenever $|E_{\text{unt},N-l}|$ depends only on the number of petals 
of the subgraph (and not the subgraph itself), we could 
infer a final formula for the number of faces. 
The method that we will use in Section 
\ref{sect:proofs} might hopefully find an extension this more general situation. Although not carried out in the present work, 
we hope that this analysis can be applied specifically for 
``periodic'' $N$-petal flowers. These ribbon graphs can be defined 
as a $\{(N_1,s_1),(N_2,s_2),(N_3,s_3),\dots,(N_q,s_q)\}$-petal flowers
(with merged sectors) with alternate sequence 
$(N_i,s_{i})=(k_1,\pm)$, $(N_{i+1},s_{i+1})=(k_2,\mp)$, $ i,i+1\in [\![1, q]\!]$, 
with fixed $k_1\geq 1$ and $k_2\geq 1$.
For any of these flowers, the class $(F,s)$ 
(hence the number of faces) can be derived. 

Let us introduce the following quantity
\beq
A(\eta,\xi,k) = \sum_{l=\eta }^\xi (-1)^{l\, k} \,, 
\eeq
where $\eta = 0,1,$ $k \geq 1$ and $\xi \geq \eta$ are
all integers.

Assume that $q=2\ell$ is even, $\ell \geq 1$, such that the signs
alternate $(s_1=-,s_2=+,\dots,s_{q-1}=-,s_q=+)$,
is such a case we have
\beq
(F,s) = 1+ A(0,\ell-1,k_1)\sum_{a=1}^{k_1}(-1)^a 
+ A(1,\ell,k_1) k_2 \mod 3\,.
\eeq
Now if $k_1$ is even, 
\beq
(F,s)=  1+ \ell k_2 \mod 3\,.
\eeq
Otherwise $k_1$ is odd and we get 
 \beq
(F,s)= 
1+ \frac{((-1)^{\ell} -1)}{2} (1+k_2) \mod 3\,.
\eeq
Let us assume that $q=2\ell +1$, $\ell \geq 1$, the sequence 
$(s_1=-,s_2=+,\dots,s_{q-1}=+,s_q=-)$ yields 
\beq
(F,s) = 1+ A(0,\ell,k_1)\sum_{a=1}^{k_1}(-1)^a 
+  A(1,\ell,k_1) k_2 \mod 3\,.
\eeq
Assuming that $k_1$ is even, then 
\beq
(F,s) = 1+ \ell k_2 \mod 3\,.
\eeq
Otherwise if $k_1$ is odd, we get 
\bea
(F,s) &=& 1+(-1) \frac{1- (-1)^{(\ell+1) k_1}}{2} 
+ (-1) \frac{1- (-1)^{k_1 \ell}}{2}k_2 \mod 3\crcr
&=&  \frac{(1- (-1)^{\ell})}{2} (1 - k_2) \mod 3\,.
\eea

Assume that $q=2\ell$ is even, $\ell \geq 1$, such that the signs
alternate $(s_1=+,s_2=-,\dots,s_{q-1}=+,s_q=-)$,
then we obtain
\beq
(F,s) = 1+A(0,\ell-1,k_1) \sum_{a=1}^{k_1}(-1)^a 
+ A(0,\ell-1,k_1) k_2 \mod 3\,.
\eeq
Now if $k_1$ is even, 
\beq
(F,s)=  1+ \ell k_2 \mod 3\,.
\eeq
But if $k_1$ is odd, one gets 
\beq
(F,s) = 1+ \frac{( (-1)^{\ell } -1)}{2} ( 1- k_2)\mod 3\,.
\eeq
The last case concerns $q=2\ell+1$, $\ell \geq 1$,
for which  $(s_1=+,s_2=-,\dots,s_{q-1}=-,s_q=+)$,
then we obtain
\beq
(F,s) = 1+A(0,\ell-1,k_1)\sum_{a=1}^{k_1}(-1)^a 
+  A(0,\ell,k_1) k_2 \mod 3\,.
\eeq
Setting $k_1$  even yields
\beq
(F,s) = 1+ (\ell +1)k_2 \mod 3\,,
\eeq
whereas assuming that $k_1$ is odd gives
\beq
(F,s) =
 \frac{1+ (-1)^{\ell}}{2} (1+k_2) \mod 3\,.
\eeq
The simplest periodic rosettes 
are defined such that 
$q=1$, $\ell =0$, and $k_1=0$ or $k_2=0$ but not
at the same time. 
They have been already treated. 
The above fomulas are again valid again and  we could 
extend these  to $k_{1,2} \geq 0$ $q\geq 1$. 
Another interesting and simple case is the one 
determined by $k_1=1=k_2$. In such a situation, 
we restrict the above results such that

$(a1):\{(s_1=-,\dots,s_q= +),q=N=2\ell\}$,
$(a2):\{(s_1=-,\dots,s_q= -),q=N=2\ell+1\}$ and $k_1$ odd: 
 \beq
(a1) \;\;
(F,s)_{-+}= (-1)^{\ell } \mod 3\, 
\qquad 
(a2)\;\;  (F,s)_{--}  =  0 \mod 3 \,,
\label{eq:f++}
\eeq
in any case, $F^{(-)}=1 + \eps(N)$;

or $(b1):\{(s_1=+,\dots,s_q= -), q=N=2\ell\}$, 
$(b2): \{(s_1=+,\dots,s_q= +),q=N=2\ell+1\}$ 
and $k_1$ odd: 
 \beq
(b1)\;\; (F,s)_{+-} = 1 \mod 3\,,\qquad 
(b2) \;\; (F,s)_{++} = 1+ (-1)^{\ell}  \mod 3\,,
\label{eq:f--}
\eeq
in any situation the number of faces is always 
$F^{(+)} = 1+ \eps_{(4)}(N)$, where $\eps_{(4)}(\alpha)=1$
if $\alpha \in 4\bbN +3$ and $\eps_{(4)}(\alpha)=0$
elsewhere.

\section{Proofs of Theorem \ref{theo:brNu},
 Corollary \ref{coro:delflo} and 
Corollary \ref{coro:terminal} }
\label{sect:proofs}

We start by proving a useful result on compositions. 

\begin{lemma}[Number of compositions with $I$ odd integers] \label{lem:comp}
Let $n,P,I$ be some integers, such that $0\le I\le P \le n$,
$0\leq d\leq D-1$. 
The number $\cC_n(P,I)$ of compositions of $n$ in $P$ integers  among which $I$ odd integers is given by
\beq
\cC_{n}(P,I) =\bin{\frac{n+I}{2}-1}{P-1}\,\bin{P}{I},
\eeq
where the brackets stand for binomial coefficients, $n-I \in 2\bbN$ (i.e. $n$ and $I$ should share the same
parity) and 
$n \geq  2P$ if $n$ is even and $n \geq 2P-1$ if $n$ is odd. 

\end{lemma}

\proof  Counting the number of compositions of the integer $n$ made of $P$ integers is well-known as $(^{n-1}_{P-1})$. 
However we now have to distinguish these compositions 
depending on how many odd integers they contain. 

The easy method to achieve this is to add $+1$ to every odd integer in the composition of $n$. 
More precisely, let us assume that the number of odd integers is $I$. The initial integer $n$ and $I$ have the same parity.
 We add $+1$ to the $I$ odd integers thus obtaining a composition of the larger integer $n+I$ in terms of solely even integers, 
i.e a straightforward composition of the integer $(n+I)/2$. Reversely, starting from a composition of the integer $(n+I)/2$, 
we multiply it by 2 and subtract $1$ to $I$ arbitrary integers among the list of $P$ integers of the composition. 
Counting the number of such possibilities, it is convenient to distinguish the cases $n$ even and odd for the sake of clarity.  
For $n=2m$ even, $I=2f$ is also even with $f$ running from 0 to the integer part of $\f P2$ (which is bounded by $m$) and we get:
\bea
\cC_n(P,I)
&=&
\Big{|}\Big\{(q_1,..,q_P)\, |\,q_l \in \bbN\setminus\{0\},\, \sum_p q_p=2m \;\textrm{and}\; \sum_p\eps(q_p)=2f\Big\}\Big{|}\crcr
&=&
\bin{m+f-1}{P-1}\,\bin{P}{2f}.
\eea
For $n=2m+1$ odd, $I=2f+1$ is odd and runs from $1$ to $P$ and we get:
\bea
\cC_n(P,I)
&=&
\Big{|}\Big\{(q_1,..,q_P) \,|\,q_l \in \bbN\setminus\{0\},\, \sum_p q_p=2m+1\;\textrm{and}\; \sum_p\eps(q_p)=2f+1\Big\}\Big{|}\crcr
&=&
\bin{m+f}{P-1}\,\bin{P}{2f+1}\,.
\eea
\qed

As a proof check of Lemma \ref{lem:comp} that summing over all values of $f$ (or equivalently $I$), 
one recovers the number of compositions $(^{n-1}_{P-1})$ of $n$ in $P$ integers. The  appendix provides more details on this development. 

The number of compositions of $n$ in $P$ parts
 which contain $I$ specific integers belonging 
to the set $D \bbN +d$, for any $d,D\in \bbN$, 
$D\geq 2$, $1\leq d\leq D-1$,
is fully addressed in the more elaborate framework 
of generating functions. We gather all these 
developments in the appendix 
from which both $\cC_n(P,I)$ and 
$\cC^t_n(P,I)$ (of Theorem \ref{theo:brNu}) can be simply deduced 
by applying $D=2,3$ and $d=1,2$, respectively, from Lemma \ref{lem:Dcomp}.

We can now address the proof of our main theorem.

\noindent \proof[Proof of Theorem \ref{theo:brNu}]

Since all the rosette graphs (hence flowers) have rank equal to 0, the BR polynomial 
of any flower $\cG_{N}$ (twisted or not) with $N$ intertwined petals reads:
\beq
R_{\cG_N}(Y,Z)=\sum_{A\subset\cG_N} Y^{n(A)}\,
Z^{ 1+n(A)- F(A)}\,,
\eeq
where $A$ runs over all spanning subgraphs of $\cG_N$. Such subgraphs can be classified as follows. 
Erasing some petals, we get subgraphs made of $P$ packets of intertwined petals, 
each packet $p$ with a number $q_p$ of petals, respectively. 
The total number of petals of the subgraph is equal to its nullity, $n=\sum_{p=1}^P q_p$, 
and takes values from 0 (the empty graph with a single vertex) to the initial nullity $N$. 
It remains to study the number of faces and in each case 
their number is different. 

Let us focus first on the untwisted flower. 
The number of faces actually depends on the parity of the number of petals of each packet.  By Theorem \ref{theo:face}, 
 defining an index $\eps(q)$ equal to 1 if $q$ is odd and 0 if $q$ is even,
\beq
\eps(q)\equiv\,\f12\big{(}1-(-1)^q\big{)},
\eeq 
the number of faces counts the number of odd numbers among the list $(q_1,..,q_P)$:
\beq
F=1+\sum_{p=1}^P \eps(q_p)\,.
\eeq
Putting all these pieces together, we get:
\bea
R_{\cG_N}(Y,Z)
&=&
\sum_{n=0}^N Y^n\, \sum_{P=0}^n \cN_{N,n}^P
\sum_{\sum_{p=1}^P q_p =n}
Z^{n-\sum_p\eps(q_p)}\label{interm} \\
&=&
\sum_{n=0}^N Y^n\, \sum_{P=0}^n \cN_{N,n}^P
\sum_{ I=0}^{P}
Z^{n-I}\,
\Big{|}\Big\{(q_1,..,q_P)\,| \sum_p q_p=n \;\textrm{and}\; \sum_p\eps(q_p)=I\Big\}\Big{|},\nonumber
\eea
where $\cN_{N,n}^P$ counts the number of occurrence that a given subgraph made of 
$P$ packets and missing $N-n$ petals can be obtained from the full graph with $N$ petals. 
Then $I$ counts the number of odd numbers in the composition $(q_1,..,q_P)$ of the integer $n$.

Easy combinatorics allows to compute $\cN_{N,n}^P$. The goal is to put back $(N-n)$ petals 
between the $P$ packets with the possibility of putting petals on the right or on left of all these packets, 
but with the constraint that at least one petal between each packet. 
Forgetting a  moment about the extremities (i.e putting 0 petals at the left or right of all packets),
 one is counting the number of compositions of $(N - n)$ in $(P-1)$ integers. 
This is given by the binomial coefficient $(^{N-n-1}_{P-2})$. Taking into account all the possibilities, we get:
\beq
\cN_{N,n}^P=\bin{N-n-1}{P-2} + 2\bin{N-n-1}{P-1}+\bin{N-n-1}{P}
\,=\,
\bin{N-n+1}{P}\,.
\eeq
We finally put the lowest and highest order monomials separate to avoid ambiguities
 in the definition of the binomial coefficients. 
The fact that $I \in [\![\eps(n),P]\!]$ comes from some 
parity constraints that should satisfy both $n$ and $I$ 
in order $\cC_n(P,I)$ to be non vanishing, see Lemma \ref{lem:comp}. 
This achieves the proof of \eqref{eq:unflow}. 

Let us discuss now the case of the twisted flower. 
The number of faces of the subgraph depends now
on the $\bbZ_3$ class corresponding to the $q_p$'s.  
From Theorem \ref{theo:face}, it is simple, to determine that
\beq
F= 1+ \sum_{p=1}^{P} \eps_{(3)}(q_p)\,, 
\eeq
where $\eps_{(3)}(q_p)=1$ if $q_p \in 3\bbN +2$
otherwise $ \eps_{(3)}(q_p)=0$. 

Thus, we obtain in similar way than before
\bea
R_{\cG^t_N} (Y,Z)&=& \sum_{n=0}^N Y^N \sum_{P=0}^n
\cN^P_{N,n} \sum_{\sum_{p}^P q_p=n} 
Z^{n - \sum_{p} \eps_{(3)}(q_p) }  \crcr
&=& \sum_{n=0}^N Y^N \sum_{P=0}^n
\cN^P_{N,n} \sum_{I =0}^{P} 
Z^{n - I }\Big{|}
\Big\{(q_1,\dots,q_P)|\, \sum_{p=1}^P q_p = n \,\text{and}\, 
\sum_{p=1}^P \eps_{(3)}(q_p) = I \Big\}\Big{|} \,,
\crcr
&&
\eea
where  $\cN^P_{N,n}$ is exactly the same as previously since 
the procedure of determining of how many subgraphs with $P$ 
packets the nullity of which is fixed to $n$ are obtained  from the initial graph remains the same. We finally get the formula \eqref{eq:twflow}
after extracting the contribution of the lowest and highest nullity subgraphs and notice that $I$ should belongs to $[3-{\rm d}_n, P]$ 
after satisfying some constraints for getting non vanishing  
$\cC^t_{n}(P,I)$, see  Lemma \ref{lem:Dcomp} in 
the appendix.  
This completes the proof of the theorem. 

\qed

We can now provide some examples of $R_{\cG_N}$ and $R_{\cG^t_N}$. 
Invoking Lemma \ref{lem:comp} in order to compute the remaining
cardinal appearing in \eqref{interm},  we get:
\bea\label{eq:rcgn}
R_{\cG_N}(Y,Z) &=&
1+Y^NZ^{N-\eps(N)}\\
& +& 
\sum_{n=1}^{N-1} Y^n\, \sum_{P\in [\![1,n]\!] \,\& \,P \leq N-n+1 } \bin{N-n+1}{P}\Bigg\{\crcr
&&
\sum_{\f {n-P}2\le k\le \min(\f n2, n-P)}
Z^{2k}\,
\bin{n-k-1}{P-1}\,\bin{P}{n-2k}\Bigg\},
\crcr
&&\nonumber
\eea
where we put the lowest and highest order monomials separate to avoid ambiguities in 
the definition of the binomial coefficients. $k$ is an integer 
and it is bounded as given in the formula.
We can check this formula explicitly for small values of $N$:
\bea
&&
R_{\cG_1}(Y,Z) = Y + 1\,,\cr\cr
&&
R_{\cG_2}(Y,Z)=Y^2Z^2+2Y+1\,, \cr\cr
&&
R_{\cG_3}(Y,Z)=Y^3Z^2+ Y^2(2Z^2+1) + 3Y+1\,, \cr\cr
&&
R_{\cG_4}(Y,Z)=Y^4Z^{4}+ 4Y^3Z^2 +  3Y^2(Z^{2} +1)
+  4Y + 1 \,,\cr\cr
&&
R_{\cG_5}(Y,Z)=Y^5 Z^{4} + Y^4(3Z^{4}+2Z^2)+
Y^3( 9 Z^{2}  +1 )
+ Y^2(4Z^{2} +6) + 5Y  + 1\,.
\eea
Meanwhile, for $\cG^t_N$, we have
\bea&&
R_{\cG^t_N}(Y,Z)= 1+Y^NZ^{N-\eps_{(3)}(N)} + 
\sum_{n=1}^{N-1} Y^n \sum_{P\in [\![1,n]\!] \,\& \,P \leq N-n+1} \bin{N-n+1}{P} \Bigg\{
\crcr
&&
\sum_{I =0}^{P} 
Z^{n - I }\Big[
\sum_{ [l \in [\![ 0, P-I ]\!] ]\,\&\, 
[ 2l \in M_{3, 3- {\rm d}_{n + I} } \text{ or } M_{3,0}]  }
\bin{P}{l,I, P- l - I} 
\bin{\frac{n +  I +2 l }{3} -1}{P-1} \Big]\Bigg\},
\label{eq:rcgn2}
\eea
where $n \geq 3P - {\rm d}_{I}$ and the last sum is only nontrivial 
for terms such that $n+I +2l \geq3 P$. Setting $N=1,2,\dots,5$, 
we obtain
\bea
&&
R_{\cG^t_1}(Y,Z) = YZ + 1\,,\cr\cr
&&
R_{\cG^t_2}(Y,Z)=
Y^2Z+2YZ+1\,, \cr\cr
&&
R_{\cG^t_3}(Y,Z)=
Y^3Z^{3} + Y^2(Z^2 + 2Z) + 3 Y Z  + 1\,, 
\cr\cr
&&
R_{\cG^t_4}(Y,Z)= Y^4 Z^{4} + 
2Y^3(Z^{3}+Z^2) + 3Y^2(Z^2 + Z) + 4 Y Z  + 1\,,
\cr\cr
&&
R_{\cG^t_5}(Y,Z)=
Y^5Z^{4} + Y^4(4Z^4 + Z^2) + Y^3(4Z^3 + 6Z^2) + Y^2(6Z^2 + 4Z) 
+ 5YZ  + 1\,.
\eea

\noindent \proof[Proof of Corollary \ref{coro:delflo}]
This is a direct consequence of the factorization 
of the BR polynomial in terms of polynomials for product
of graphs.   
We recall first that for a product graph $\cG_1 \cdot \cG_2$ of two disjoint ribbon graphs $\cG_1$ and $\cG_2$ glued along
one of their vertex \cite{bollo}, the BR polynomial 
is given by 
\beq
R_{\cG_1 \cdot \cG_2}  =  R_{\cG_1} \cdot R_{\cG_2}\,. 
\label{eq:prod}
\eeq
Consider now a graph $\cG$ having a trivial $N$-petal flower $G_N$
subgraph, then it is easy to factor $\cG$ as $G_N \cdot \cG_{\cG - G_N}$. 
The statement is therefore obvious from \eqref{eq:prod}. 
\qed

\noindent \proof[Proof of Corollary \ref{coro:terminal}]
Using Corollary \ref{coro:delflo} and  previous results
on terminal forms as bridges and trivial self-loops, 
the result becomes immediate. 
\qed

\medskip
 
\noindent {\bf On a $(1,1)$-periodic $N$-petal flower.}
Consider the $(k_1,k_2)$-periodic $N$-petal 
flower  (according to our discussion in Section \ref{sect:gene}).
The problem of finding an explicit expression for 
the BR polynomial associated with this graph becomes 
more involved and should reduce again to a counting of 
specific compositions.  

For a matter of simplicity, let us discuss the case
$(k_1 = 1, k_2= 1)$-periodic $N$-petal flower $\cG'_N$.  
We again consider an expansion of the polynomial 
in terms of the nullity. At fixed nullity $n$, we consider 
$P$ packet subgraphs $(q_1,\dots,q_P)$. 
According to \eqref{eq:f++} and \eqref{eq:f--}, 
the number of faces in each packet  depends on the quality 
of the its last petal. It is clear that we can decompose 
the list $(q_1,..,q_P)$ in $(q_{i_1},..,q_{i_{l_1}})$
for which the last petal  of each $i_p$ packet is untwisted,
and $(q_{j_1},..,q_{j_{l_2}})$ for which 
the last petal  of each $j_p$ packet is twisted. 
Hence $l_1+l_2 = P$. 

The number of faces in the subgraph corresponding to 
the list $(q_1,..,q_P)$ is simply given by 
\beq
F=1+\sum_{p=1}^{l_1} \eps(q_{i_p})
+\sum_{p=1}^{l_2} \eps_{(4)}(q_{j_p})\,.
\eeq
Arguing as above in the above proof of Theorem \ref{theo:brNu}, 
we see that there is two layers of difficulty when one 
asks for a closed formula for a polynomial for $\cG'_N$. 
For a subgraph corresponding  to 
a composition $(q_1,..,q_P)$, we need to know which of
the packets start by a twisted petal (and so which do not).
This leads to consider more general ``signed'' composition 
$(q_1^{\varepsilon_{1}}, q_2^{\varepsilon_{2}},
\dots,q_P^{\varepsilon_{P}})$, where $\varepsilon_{i} =\pm$, determines
if the starting petal in the packet $i$ is twisted or not, respectively. 
Having find a way to encode these  
signed compositions for any subgraph, and summing on 
all of these subsets of compositions, 
then the remaining task is to count again among these compositions those containing $I$ odd integers and $J$ integers in $4\bbN + 3$.

\section{Recurrence relations}
\label{sect:recrel}
In fact, there is an alternative way to determine the above 
BR polynomial of the $N$-petal flower worthwhile to be
discussed.  

Consider a $N$-petal flower $\cG$ and its spanning subgraphs.
We can assign an index to each petal, $\{e_i\}_{i=1,\dots, N}$
and we can, for instance, label them from the top to the bottom
(using the cyclic order). Then, choose $e_1$ (or $e_N$, without loss of generality). 
Then the spanning subgraphs of $\cG$ can be divided
into two  sets: a set $S_{e_1}$ which contains $e_1$ 
and another $S_{\check{e}_1}$ which does not. It is immediate that 
the polynomial for $\cG$ can be written 
\beq
R_{N}  = R_{N-1} + R'
\eeq   
where $R_{N-1}$ is a BR polynomial for a $(N-1)$-petal flower 
built from $S_{\check{e}_1}$ or spanning subgraphs 
made with $\{e_2,e_3,\dots, e_N\}$ and $R'$ is a sum of all contributions 
coming from subgraphs containing $e_1$.   
Among these subgraphs containing $e_1$, we can
specify subgraphs which do not contain $e_2$  
and those which do. A moment of thought leads one
to the result
\beq
R' = c(Y,Z) R_{N-2} + R''
\eeq
where $c(Y,Z)$ depends on the type of the petal $e_1$;
if $e_1$ is untwisted then $c(Y,Z) = Y$; 
if $e_1$ is twisted then $c(Y,Z) = YZ$. 
One iterates the procedure until there is no petal 
left. Finally, adding the empty graph and the total 
graph contribution, we get a recurrence
relation that $R_N$ should satisfy:
\beq
R_{N}(Y,Z)  = \sum_{n=0}^{N-1} Y^n Z^{n- \eps_{(\alpha)}(n)}
R_{N-n-1}(Y,Z) + Y^N Z^{N- \eps_{(\alpha)}(N)}\,,
\eeq 
for all $N\in \bbN$, $N \geq 1$, and where $\alpha=2,3$ given
the particular type of flower we are dealing with: 
$\eps_{(2)}(q) = (1-(-1)^q)/2$ if we have a untwisted
flower and $\eps_{(3)}(q)=1$ if $q \in 3\bbN + 2$
otherwise $\eps_{(3)}(q)=0$, for the twisted case. 
  
Notice that such a recurrence relation is difficult to solve
for arbitrary $N$  by ordinary methods.  
However, given $\alpha=2,3$, the above techniques
show that explicit solutions for such recurrence equations
are affordable.

\begin{center}
{\bf Acknowledgements}
\end{center}

{\footnotesize Discussions with Vincent Rivasseau are gratefully acknowledged. RCA thanks the Perimeter Institute for its hospitality.
Research at Perimeter Institute is supported by the Government of Canada through Industry
Canada and by the Province of Ontario through the Ministry of Research and Innovation.
This work is partially supported by the Abdus Salam International Centre for Theoretical Physics (ICTP, Trieste, Italy) through the Office of External Activities (OEA)-Prj-15. The ICMPA is also in partnership with the Daniel Iagolnitzer Foundation (DIF), France.}

\appendix
\section*{Appendix: On compositions}

\renewcommand{\theequation}{A.\arabic{equation}}
\setcounter{equation}{0}

We address in this appendix  important remarks on the number of compositions
and number of specific compositions containing a certain number of
particular integers. These numbers of compositions have been extensively used in the text.

 The number of compositions of the integer $n$ in $P$ integers is the coefficient of $t^n$ of the following generating function:
\bea
\label{eq:compo}
\left(\sum_{k=1}^{\infty}t^k\right)^P
 &=& \sum_{n=1}^{\infty}\Big(\sum_{k_1 +k_2 +\dots +k_P= n} 1\Big)
t^n 
\crcr
& =& t^P\f1{(1-t)^P}
= t^P\sum_{k=0}^\infty t^k \bin{P+k-1}{k}
=\sum_{n=P}^\infty \bin{n-1}{P-1}  t^n\, .
\eea
In order to compute the number of compositions of $n$ in $P$ 
parts containing a number $I$ of odd integers (or of $P-I$ even integers), we further decompose the above expression as
\bea
\left(\sum_{k=1}^{\infty}t^k\right)^P
 = \sum_{n=1}^{\infty}\sum_{I=0}^P
\Big(\sum_{k_1 +k_2 +\dots +k_P= n \,\&\, I_2(k_1,k_2,\dots,k_P)=I} 1\Big)
t^n \,,
\eea
where $(k_1,k_2,\dots,k_P)$ is a composition of $n$, 
$I_2(k_1,k_2,\dots,k_P)$ is a function calculating the number
of odd integers in that composition (respectively, calculating the number of even integers 
in the composition; in the above sum the condition should change as $I_2(k_1,k_2,\dots,k_P)=P-I$).
 Next, we perform the following transformation:
\bea
\left(\sum_{k=1}^{\infty}t^k\right)^P
 &=& \sum_{I=0}^P \bin{P}{I}
\left(\sum_{k\in 2\mathbb{N} +1}^{\infty}t^k\right)^I
\left(\sum_{k\in 2\mathbb{N} +2}t^k\right)^{P-I}
=
 \sum_{I=0}^P \bin{P}{I} t^{-I}
\left(\sum_{k\in 2\mathbb{N} +2}^{\infty}t^k\right)^P\crcr
&=& \sum_{I=0}^P \bin{P}{I}
\sum_{k=0}^\infty \bin{P+k -1}{P-1} t^{2(k+P) - I} \crcr
&=&
\sum_{n\in 2\mathbb{N}\,\&\, n \geq 2P}^\infty 
\sum_{I'=0}^{[\f P2]} \bin{P}{2I'}
 \bin{ \f {n+2I'}{2} -1 }{P-1} t^{n}
 \crcr
&& 
+ \sum_{n\in 2\mathbb{N}+1\,\&\, n \geq 2P-1}^\infty
\sum_{I'=0}^{[\f P2]} \bin{P}{2I'+1}
 \bin{\f {n+2I'+1}{2} -1}{P-1} t^{n}\,,
\eea
where we use \eqref{eq:compo} as an intermediate step
to compute the sum over even integers. 
Then, we get
\bea
&&
n \in 2\mathbb{N} \;\,\text{and}\;\, n\geq 2P\,,
\qquad 
\sum_{I=0\,\&\, I\, \text{even}}^{P} \bin{P}{I}
 \bin{ \f {n+I}{2} -1 }{P-1}  = \bin{n-1}{P-1} ;\crcr
&&
n \in 2\mathbb{N}+1 \;\,\text{and}\;\, n\geq 2P-1\,,
\qquad 
\sum_{I=1\,\&\, I\, \text{odd}}^{P} \bin{P}{I}
 \bin{ \f {n+I}{2} -1 }{P-1}  = \bin{n-1}{P-1} .
\eea
The number of compositions containing $I$ odd 
integers is simply given by the summand  of the above
$\cC_n(P,I)=(^P_I)\big(^{(n+I)/2 -1 }_{P-1}\big)$. 
 
It is far advantageous to introduce a generalized version of Lemma \ref{lem:comp} valid in any situation. 

\begin{lemma}[Number of compositions with $I$ integers belonging to 
$D\bbN +d$] \label{lem:Dcomp}
Let $n,P,I,D,d$ be some integers, such that $I\le P \le n$,
$1\leq d\leq D-1$, $D \geq 2$.  
We denote ${\rm d}_q$ the remainder 
of the Euclidean division of $q\in \bbN$ by $D$.
 The number $\cC^{D,d}_{n,P,I}$ of compositions of $n$ in $P$ parts among which $I$ integers belonging to $D \bbN +d$ is given by 
\bea
&&
\cC^{D,d}_{n,P,I}=  
 \sum_{[\sum_{\alpha=1 \,\&\,\alpha \neq d}^{D} l_\alpha =P-I ]\,\&\, [\sum_{\alpha=1 \,\&\,\alpha \neq d}^{D-1} (D-\alpha) l_\alpha 
\in M_{D, D - {\rm d}_{n+(D-d)I} } \text{ or } M_{D,0}]}
\Bigg\{\crcr
&&
\bin{P}{l_1,l_2,\dots, l_{d-1},I,l_{d+1} \dots,l_{D}} 
\bin{\frac{n + (D-d) I +\sum_{\alpha=1 \& \alpha \neq d}^{D-1} (D-\alpha) l_\alpha }{D} -1}{P-1}\Bigg\},
\eea 
where the first and  the second bracket denote 
the multinomial and binomial coefficients, respectively.  
\end{lemma}

\proof

Let $D,d\in \bbN$, such that $d\le D$.  Consider 
the set 
\beq
M_{D,d} = \{d, D+d, 2D+d, 3D+d, \dots, kD +d, \dots\}
 = D \bbN + d  \subset \bbN\,.
\eeq
Remark that $M_{0,0}=\{0\}$,
 $M_{1,0}=\bbN$, $M_{2,1}$ is the set of odd positive integers 
and $ M_{2,0}$ is the set of  even positive integers including 0. 
In general, $M_{D,0}$ is the set of  multiples of $D$
and $\bbN =\cup_{d=0}^{D-1} M_{D,d}$ if $D \geq 1$.

Let us define the function $I_{D,d}$ on the composition $(k_1,k_2, \dots, k_P)$ of some integer $n$ in $P$ integer parts, such that 
$I_{D,d}(k_1,k_2, \dots, k_P)$ counts the numbers of elements 
$k_i$ of the composition which belong to $M_{D,d}$. Note that 
$I_{2,1}$ counts the number of odd integers in the composition,
hence the above notation $I_2= I_{2,1}$. 

Expanding the generating function as 
\bea
\left(\sum_{k=1}^{\infty}t^k\right)^P
 = \sum_{n=1}^{\infty}\sum_{I=0}^P
\Big(\sum_{k_1 +k_2 +\dots +k_P= n\,\&\, I_{D,d}(k_1,k_2,\dots,k_P)=I} 1\Big)
t^n \,,
\eea
for $D\geq 2$, $1\leq d \leq D-1$, we aim at computing the coefficient of $t^n$. We use the same technique as before and find,
for all $D\geq 2$, 
\bea
&&
\left(\sum_{k=1}^{\infty}t^k\right)^P
 =
\left(\sum_{k\in M_{D,D}}^{\infty}t^k + 
\sum_{k\in M_{D,D-1}}^{\infty}t^k + 
\dots +
\sum_{k\in M_{D,2}}^{\infty}t^k +
\sum_{k\in M_{D,1}}^{\infty}t^k 
\right)^P
\crcr
&&=\sum_{l_1+l_2+\dots+l_{D} =P}
\bin{P}{l_1,l_2, \dots,l_{D-1},l_{D}}
 \crcr
&&
\left(\sum_{k \in M_{D,D}}^{\infty}t^k\right)^{l_D}
\left(\sum_{ k \in M_{D,D-1}}^{\infty}t^k\right)^{l_{D-1}}
\dots
\left(\sum_{ k \in M_{D,2}}t^k\right)^{l_{2}}
\left(\sum_{k \in M_{D,1}}^{\infty}t^k\right)^{l_1},
\eea
where the first bracket $(^P_{l_1,\dots,l_{D}})$ stands for a multinomial coefficient. Computing further, we obtain
\beq
\left(\sum_{k=1}^{\infty}t^k\right)^P
=
\sum_{l_1+l_2+\dots+l_{D} =P}
\bin{P}{l_1,l_2, \dots,l_{D-1},l_{D}}\frac{1}{t^{l_{D-1}+ 2l_{D-2}+ \dots
+ (D-1)l_1}}
\left(\sum_{k\in M_{D,D}}^{\infty}t^k\right)^{P},
\eeq
and using again \eqref{eq:compo}, we write
\beq
\left(\sum_{k=1}^{\infty}t^k\right)^P=
\sum_{l_1+l_2+\dots+l_{D} =P}
\bin{P}{l_1,l_2, \dots,l_{D-1},l_{D}}
\frac{1}{
t^{\sum_{\alpha=1}^{D-1}
(D-\alpha) l_{\alpha} } }\sum_{k=0}^{\infty}
\bin{P+k -1}{k} t^{D(k+P)} .
\eeq
We fix $1\leq d \leq D-1$, $l_d=I$ and want to perform a change in variable as 
\beq
k \to n = D(k+P) -(D-d) I-\sum_{\alpha=1\,\&\, \alpha \neq d}^{D-1} (D-\alpha) l_\alpha\,,
\label{eq:kn}
\eeq
this requires that 
$n + (D-d) I + \sum_{\alpha=1\,\&\, \alpha \neq d}^{D-1}(D-\alpha) l_\alpha$ should be a multiple of $D$. Given $q \in M_{D,s}$, we denote ${\rm d}_{q} =s$, therefore \eqref{eq:kn} implies that 

- ${\rm d}_{ \sum_{\alpha=1\,\&\, \alpha \neq i}^{D-1}(D-\alpha) l_\alpha} =  D - {\rm d}_{(n+ (D-d) I)} < D $ 
which holds for ${\rm d}_{(n+(D-d)I)}$ non vanishing;

-  ${\rm d}_{ \sum_{\alpha=1\,\&\, \alpha \neq d}^{D-1}(D-\alpha) l_\alpha} = 0$
for ${\rm d}_{(n+ (D-d) I)}$ vanishing. 

One notices that, given $n,D,d$ and $I$,  these two circumstances exclude
each other. We write these conditions 
$[\sum_{\alpha=1 \,\&\,\alpha \neq d}^{D-1} (D-\alpha) l_\alpha 
\in M_{D, D - {\rm d}_{(n+(D-d)I)} } \text{ or } M_{D,0}]$
but one pays attention that, given $n,D$ and $d$, 
one has to choose between these. 

Now, the change the order of the summations and
introduce some constraints in the $l_\alpha$, 
for a given $1\leq d \leq D-1$,
\bea
&&
\left(\sum_{k=1}^{\infty}t^k\right)^P = \crcr
& &=\sum_{ k=0}^\infty\;\,
\,\sum_{I=0}^{P} \,
 \sum_{\sum_{\alpha=1 \,\&\,\alpha \neq d}^{D} l_\alpha =P-I }
\Bigg\{ \crcr
&&
\bin{P}{l_1,l_2,\dots, l_{d-1},I,l_{d+1} \dots,l_{D}} 
\bin{P+k -1}{P-1} t^{D(k+P) -(D-d) I-\sum_{\alpha=1\,\&\,\alpha \neq d}^{D-1} (D-\alpha) l_\alpha} \Bigg\} \,.\crcr
&&
\label{intermD}
\eea
Hence, changing now variable $k \to n$, 
$n \geq DP - {\rm d}_{(D-d)I}$, the coefficient
of $t^n$ in the above sum is given by 
\bea
&&
\sum_{I=0}^P 
 \sum_{[\sum_{\alpha=1 \,\&\,\alpha \neq d}^{D} l_\alpha =P-I ]\,\&\, [\sum_{\alpha=1 \,\&\,\alpha \neq d}^{D-1} (D-\alpha) l_\alpha 
\in M_{D, D -  {\rm d}_{(n + (D-d)I) }| } \text{ or } M_{D,0}]}
\Bigg\{\\
&&
\bin{P}{l_1,l_2,\dots, l_{d-1},I,l_{d+1} \dots,l_{D}} 
\bin{\frac{n + (D-d) I +\sum_{\alpha=1 \& \alpha \neq d}^{D-1} (D-\alpha) l_\alpha }{D} -1}{P-1}  \Bigg\} = \bin{n-1}{P-1}\,.
\nonumber
\eea
We also deduce that the cardinal of 
\bea
  \Big\{(q_1,\dots,q_P) | \; q_l \in \bbN\;, \sum_{l=1}^P q_l = n 
\;\text{ and }\; I_{D,d}(q_1,\dots,q_P)= I\Big\}\crcr
\eea
is nothing but
\bea
&&
\cC^{D,d}_{n,P,I}  = \sum_{[\sum_{\alpha=1 \,\&\,\alpha \neq d}^{D} l_\alpha =P-I ]\,\&\, [\sum_{\alpha=1 \,\&\,\alpha \neq d}^{D-1} (D-\alpha) l_\alpha 
\in M_{D, D -  {\rm d}_{(n + (D-d)I)} } \text{ or } M_{D,0}]}
\Bigg\{ \crcr
&&
\bin{P}{l_1,l_2,\dots, l_{d-1},I,l_{d+1} \dots,l_{D}} 
\bin{\frac{n + (D-d) I +\sum_{\alpha=1 \& \alpha \neq d}^{D-1} (D-\alpha) l_\alpha }{D} -1}{P-1} \Bigg\}
\eea
for $n \geq DP - {\rm d}_{(D-d)I}$. 

\qed

We emphasize that it could happen that 
 the binomial  coefficient vanishes.   
The additional condition for the coefficient  to be 
non vanishing reads off as
\beq
n + (D-d) I + \sum_{\alpha=1\,\&\, \alpha \neq d}^{D-1}(D-\alpha) l_\alpha \geq DP \,.
\eeq
 
In particular,  computing 

- the $N$-untwisted-petal  flower
we need to evaluate $\cC^{2,1}_{n,P,I}$ (and denoted by $\cC_{n}(P,I)$ in 
Theorem \ref{theo:brNu}) that we can deduce from 
the above formula as, given ${\rm d}_n=d$,
\bea
\cC^{2,1}_{n;P;I}   =  \bin{P}{I,P-I}
\bin{\frac{n + I  }{2} -1}{P-1}=\bin{P}{I}
\bin{\frac{n + I  }{2} -1}{P-1}\,,
\eea
where $n$ and $I$ should have the same parity 
and $n \geq 2P - {\rm d}_{I}$;

- the $N$-twisted-petal flower
the relevant cardinal is $\cC^{3,2}_{n,P,I}$ 
(and denoted by $\cC^t_{n}(P,I)$ in 
Theorem \ref{theo:brNu}) that we can deduce as, given ${\rm d}_n=d$,
\bea
\cC^{3,2}_{n,P,I} &=& \sum_{[ l_1+l_3=P-I ]\,\&\, [2l_1
\in M_{3, 3  - {\rm d}_{(n+I)} } \text{ or } M_{3,0}]}
\bin{P}{l_1,I,l_3 } 
\bin{\frac{n + I +2l_1 }{3} -1}{P-1} \crcr
&=&
 \sum_{ [l \in [\![ 0, P-I ]\!] ]\,\&\, 
[ 2l \in M_{3, 3 - {\rm d}_{(n + I)} } \text{ or } M_{3,0}]  }
\bin{P}{l,I, P- l - I} 
\bin{\frac{n +  I +2 l }{3} -1}{P-1},
\eea
where $n \geq 3P - {\rm d}_{I}$ and $n+I+2l\geq 3P$.

\vspace{0.5cm}

\end{document}